\documentclass{article}

\usepackage{amssymb}
\usepackage[cp1251]{inputenc}
\usepackage[english]{babel}
\usepackage{amsmath}
\usepackage{color}
\usepackage{graphics}
\usepackage{epsfig}
\usepackage[all,knot]{xy}

\def\mapr#1{\smash{\mathop{\buildrel{#1}\over\longrightarrow}}}

\def\mapd#1{\Big\downarrow\rlap{$\vcenter{\hbox{$#1$}}$}}
\def\mapu#1{\Big\uparrow\rlap{$\vcenter{\hbox{$#1$}}$}}

\newtheorem{theorem}{Theorem}
\newtheorem{lemma}{Lemma}
\newtheorem{definition}{Definition}

\newtheorem{conjecture}{Conjecture}

\def\proof{{\bf Proof.}}%\nobreak\noindent}
\def\C{{\bf C}}

\def\K{{\bf K}}

\def\Q{{\bf Q}}
\def\R{{\bf R}}
\def\SS{{\bf S}}

\def\U{{\bf U}}

\def\Z{{\bf Z}}

\def\0{{\bf 0}}
\def\1{{\bf 1}}

\def\cA{{\cal A}}

\def\cH{{\cal H}}

\def\cK{{\cal K}}

\def\cQ{{\cal Q}}

\def\cV{{\cal V}}
\def\cW{{\cal W}}

\def\rF{{\frak F}}

\def\rH{{\frak H}}

\def\bydef{\stackrel{def}{=}}

\def\Hom{{\hbox{\bf Hom}\;}}

\def\id{{\hbox{\bf Id}}}

\def\mat{\mbox{\bf Mat}}

\def\sign{\hbox{\bf sign}\;}

\def\formulacolor{\color{black}}
\def\blf{\formulacolor $}
\def\elf{$ \color{black}}

\def\bdf{\formulacolor$$}
\def\edf{$$\color{black}}

\def\beq#1{\formulacolor\begin{equation}\label{#1}}
\def\eeq{\end{equation}\color{black}}

\def\feq#1#2{\beq{#1}#2\eeq}
\def\ff#1{\bdf #1 \edf}
\def\f#1{\blf #1\elf}

\def\mat{\formulacolor $}
\def\tam{$ \color{black}}

\newtheorem{corollary}{Corollary}
\def\({\left(}
\def\){\right)}
\def\<{\langle}
\def\>{\rangle}
\def\ch{{\hbox{ch}}}

\title{Bordisms of manifolds with proper action
of a discrete group: signatures and
descriptions of $G$-bundles}

\author{Mishchenko, A.S. \\(Moscow State University)
\thanks{Partly supported by RFBR 11-01-00057-a,
10-01-92601-KO\_a, 11-01-90413-Ukr\_f\_a and State program RNP 2.1.1.5055}
\and Morales Mel\'endez, Quitzeh \\(National Autonomous University of Mexico)
\thanks{Partly supported by CONACyT Project No.98697}}
\date{}

\begin{document}
\maketitle

\begin{abstract}
In this work  the equivariant
signature of a manifold with proper action of a discrete
group is defined as an invariant of equivariant bordisms. It is shown
that the computation of this signature can be reduced to
its computation on fixed points sets equipped with their
tubular neighborhoods. It is given a description of the equivariant
vector bundles with action of a discrete group $G$ for the
case when the action over the base is proper quasi-free,
i.e. the stationary subgroup of any point is finite.
The description is given in terms of some
classifying space.
\end{abstract}

\section{Proper action of discrete group}
%\subsection{Определения}
%In this section we give some definitions,
%clarify some facts and notions,
%and introduce notation.

\begin{definition}
Let $G$ be a discrete (countable) group,  $M$ be a smooth
orientable manifold. The action of the group $G$ on $M$ $$
G\times M \mapr{} M
$$
is called proper, if for any point ${x\in M}$ the stationary group
 ${G_{x}\subset G}$,\linebreak
${G_{x}\bydef \{g\in G: gx=x\}}$ is finite. The action is called cocompact
if the quotient space ${M/G}$ is compact.
\end{definition}

As a more accurate definition, we have the following.
\begin{definition}
The action of the group ${G}$ on ${M}$ is proper if the mapping
$${G\times M\mapr{} M\times M, \quad (g,x)\longmapsto (gx,x),
\quad g\in G, x\in M,
}$$
is proper.
\end{definition}

%\subsection{Симплициальные и гладкие действия}

The action is called smooth if each element
${g\in G}$ acts on ${M}$ by a diffeomorphism.

Let ${H < G}$ be a finite subgroup.
Denote by ${M^{H}}$ the subset of fixed points
$${M^{H}\bydef \{x\in M: \forall g\in H\; gx=x\}.}$$

The subset ${M^{H}}$ is a smooth submanifold, see for example
\cite{Paw-02e}

Let ${\cH(G)}$ be the family of all (finite) subgroups
${H\subset G}$ such that ${M^{H}\not = \emptyset}$.
The family ${\cH(G)}$ forms a category with respect to
embedding of subgroups.
Let  ${\cH(H)\subset \cH(G)}$ be the subfamily of subgroups
 ${H'\in \cH(G)}$ such that  ${H'\supset H, \quad H'\neq H}$.

If ${H_{1}\subset H_{2}\subset G}$, then
${
M^{H_{2}}\subset M^{H_{1}}.
}$
The family
$$
V^{H}\bydef M^{H}\setminus \(\bigcup\limits_{H'\in \cH(H)}M^{H'}\).
$$
forms a smooth stratification on the manifold ${M}$.

\subsection{Simplicial proper actions}

S.Illman and T.Korppi proved that each smooth proper action is simplicial
with respect to a simplicial structure on the manifold $M$
(\cite{Kor-05en,Ill-00e}).

For a manifold with proper simplicial action one has a diagram:

\beq{1a}
\xymatrix{
C_{0}&\ar[l]_{d_{1}}C_{1}&\ar[l]_{d_{2}}\cdots&\ar[l]_{d_{n-1}}C_{n-1}&\ar[l]_{d_{n}}C_{n} \\
C^{n}\ar[u]_{D_{0}}&C^{n-1}\ar[l]^{\delta_{n}}\ar[u]_{D_{1}}&
\cdots\ar[l]^{\delta_{n-1}}&C^{1}\ar[l]^{\delta_{2}}\ar[u]_{D_{n-1}}&
C^{0}\ar[l]^{\delta_{1}}\ar[u]_{D_{n}} \\
 }
\eeq
where $${C_{k}\bydef C_{k}(M)},\quad 0\leq k\leq n,$$
and $${C^{k}\bydef \Hom_{0}(C_{k}(M),\cK),\quad 0\leq k\leq n},$$
are chain and cochain groups correspondingly.

For the diagram
\ref{1a}
the following properties hold

\begin{enumerate}
\item \mat{d_{k-1}d_{k}=0}\tam,
\item \mat{d_{k}D_{k}+(-1)^{k+1}D_{k-1}d^{*}_{n-k+1}=0}\tam,
\item \mat{D_{k}=(-1)^{k(n-k)}D^{*}_{n-k}}\tam.
\item \mat{H(D_{k}): H(C^{n-k})\mapr{}H(C_{k})}\tam is isomorphism in homology.
\end{enumerate}

The diagram \ref{1a} with properties (1--4) is called algebraic Poincare complex (APC(M,G) ).

Unlike of non simply connected manifolds, in the case of proper action
all modules  $${C_{k}\bydef C_{k}(M)}$$
and $${C^{k}\bydef \Hom_{0}(C_{k}(M),\cK)}$$
are neither free nor projective but they are finitely generated
over the group ring ${\cA=\cK[G]}$ in the case ${\cK\approx \Z}$.

A crucial point is that, in the case when ${\cK}$ is a field of zero characteristic
(for example if ${\cK\approx \Q,\R,\C}$), one has

\begin{theorem}If the action of the group $G$ is proper (cocompact) then
all modules ${C_{k}}$ over a field ${\cK}$ of zero characteristic
are projective finitely generated modules and
$${
C^{k}\approx \Hom_{\cK}(C_{k},\cK).
}$$
\end{theorem}

\subsection{Noncommutative signature for manifolds with proper action}

Let ${\K_{p}(\cA)}$ be the groups of Hermitian
${\K}$--theory, based on finitely generated projective
modules. There is a natural map
$${
\K(\cA) \mapr{} \K_{p}(\cA).
}$$

Hence one can define a noncommutative signature for manifolds with
proper action of a discrete group similar to classical case:
$${\sign M \bydef \sign(ACP(M,G))\in \K_{p}(\cA).}$$

One can reduce the algebra $\cA=\cK[G]$ to its completion $\bar{\cA}=C^*[G]$. Then
one has

\begin{theorem}
If  ${\cA}$ is a ${C^{*}}$--algebra then
${\K(\cA) \approx \K_{p}(\cA).}$

\end{theorem}

%\subsection{Бордизмы многообразий с собственным действием}

\begin{definition}
Consider an orientable manifold ${X}$ with proper action of the group ${G}$.
Say that the manifold ${X}$ is bordant to zero if there exists an orientable
manifold with boundary  ${W}$ with proper action of the group ${G}$ such that
$${
\partial W = X.
}$$
Let ${{\ }_{p}\Omega_{n}^{G}}$ be the family of all bordant classes of manifolds with
proper action of the group ${G}$.
\end{definition}

%\begin{block}{Теорема}
The noncommutative signature ${\sign(ACP(M,G))}$ satisfies the properties:

\begin{itemize}
\item It is homotopy invariant;
\item It is an invariant of bordisms.
\end{itemize}
So one has a map
$${
\sign: {\ }_{p}\Omega_{n}^{G} \mapr{} \K_{p}(\cA).}$$
%\end{block}

Let ${{\ }_{f}\Omega_{n}^{G}}$ be the similar bordism group
for free cocompact action of the group ${G}$.
One has the commutative diagram
$${\begin{array}{ccc}
      {}_{p}\Omega_{n}^{G} & \mapr{\sign} & \K_{p}(\cA) \\
      \cup &  & \cup \\
      {}_{f}\Omega_{n}^{G} & \mapr{\sign} & \K(\cA)
    \end{array}
}$$

Let ${BG}$ be classifying space, that is the Eilenberg-MacLane complex
${BG=\K(G,1)}$. Then each manifold
${M}$ with free cocompact action of the group ${G}$ induces a continuous mapping
$${\varphi_{M}: M/G\mapr{} BG},$$ such that the diagram
$${\begin{array}{ccc}
      M & \mapr{\widetilde\varphi_{M}} & EG \\
      \mapd{} &  & \mapd{} \\
      M/G & \mapr{\varphi_{M}} & BG
    \end{array}
}$$
is commutative.

Then one has
$${ {}_{f}\Omega_{n}^{G}\approx
\Omega_{n}^{G}(BG) }$$
that is the signature induces a mapping
$${
\Omega_{n}^{G}(BG)  \mapr{\sign}  \K(\cA),
}$$
and can be expressed in the terms of characteristic classes,
which is known as the Hirzebruch formula.

\subsection{The Hirzebruch formula}
%In this section we state some facts and conjectures.

The formula states that for any orientable compact manifold ${M}$
$${\sign M =\<L(M),[M] \> ,}$$
where ${L(M)}$ is the Hirzebruch characteristic class
as a polynomial of the Pontryagin classes.
If  ${M}$ is a manifold with free cocompact action of the group
${G}$, then one has a similar formula:
$${
\sign M =\<L(M/G)\ch_{\cA}(\xi_{\cA}),[M/G] \> \in \K(\cA)\otimes \cQ,
}$$
where ${\cA=C^{*}[G]}$ is the group ${C^{*}}$--algebra of the group ${G}$.

%\subsection{Классифицирующее пространство для собственного действия}

\begin{theorem}[P.Baum, A.Connes, and N.Higson.]
For a discrete group ${G}$ there is an universal space
${{}_{p}EG}$ with proper action of the group
${G}$ such that for any space ${X}$ with proper action
of the group ${G}$ there is an equivariant continuous
mapping
$${ f: X \rightarrow {}_{p}EG}$$
unique up to equivariant homotopies \cite{BaCoHi-94e}.
\end{theorem}

This means that
$${{}_{p}\Omega^{G}_{n}\approx{}_{p}\Omega^{G}_{n}({}_{p}EG).}$$

Each manifold ${M}$ with proper action of the group
${G}$ induces a stratified manifold with singularities
${M'=M/G}$. In other words one has a natural mapping
$${{}_{p}\Omega^{G}_{n}\mapr{/}{}_{s}\Omega_{n}({}_{p}EG/G).}$$

\section{Problem}
%%\subsection{Open problems}

\begin{itemize}
\item Find a Hirzebruch type formula for the noncommutative signature
for manifolds with proper action.
\end{itemize}

Let $E^{f}_{G}$ be a contractible space with free action of the group $G$.
Then, to each equivariant vector bundle $\xi$ over $M$
corresponds the bundle $pr_{1}^{*}(\xi)$ over the space
$M\times E^{f}_{G}$. Since the action of the group $G$ is free on the product
$M\times E^{f}_{G}$, then there is a vector bundle $\eta$ over the quotient space
$(M\times E^{f}_{G})/G$ such that $q^{*}(\eta)=pr_{1}^{*}(\xi)$.
Notice that the inverse image ${pr'}_{1}^{-1}(x)$ has finite homology for any point
$x\in M'=M/G$. Therefore there should be a vector bundle $\zeta$ over $M'$ such that
$(pr'_{1})^{*}(\zeta) = N\eta$ for a proper choice of the integer $N$.
This means that the
diagram
$${
\xymatrix{
{\cV ect}^{G}(M)\ar[r]^{pr_{1}^{*}}\ar@{.>}[d]& \cV ect^{G}\( M\times E^{f}_{G}\)\\
\cV ect(M'=M/G)\ar[r]_{{pr'}^{*}_{1}}^{\approx^{\otimes\cQ}} & \cV ect\(\(M\times E^{f}_{G}\)/G\)\ar[u]_{q^{*}}^{\approx}
}
}$$
is commutative and that $(pr'_{1})^{*}$ is an isomorphism up to integer
multiplication.

$${
\xymatrix{
M\ar[d]_{p}& M\times E^{f}_{G}\ar[d]^{q}\ar[l]_{pr_{1}}\\
M'=M/G & \(M\times E^{f}_{G}\)/G\ar[l]^{pr'_{1}}
}
}$$

Let us check that cohomology of the inverse image of the point, ${pr'}_{1}^{-1}(x)$
is finite. For ${x\in M^{H}}$
one has
$${
{pr'}_{1}^{-1}(x)=\(G/H\times E^{f}_{G}\)/G
}$$

$${
\xymatrix{
G/H\times E^{f}_{G}\ar[r]^-{H}&\(G/H\times E^{f}_{G}\)/H\ar[r]^{G/H}\ar[d]^{=}
&\(G/H\times E^{f}_{G}\)/G\ar[d]^{\approx}\\
&G/H\times(E^{f}_{G})/H\ar[r]_{G/H}&(E^{f}_{G})/G
}}$$

Hence
$${
H^{*}\{{pr'}_{1}^{-1}(x),x_{0};\Q\}=0
}$$
and
$$
\xymatrixcolsep{3pc}\xymatrix{
H^{*}(M/G;\Q)\ar[r]^-\approx & H^{*}\(\(M\times E^{f}_{G}\)/G;\Q\)
},$$
$$
\xymatrix{
\cV ect(M'=M/G)\otimes\cQ\ar[r]^-\approx &
\cV ect\(\(M\times E^{f}_{G}\)/G\)\otimes\cQ.
}$$

So
$${
\xymatrix{
{\cV ect}^{G}(M)\ar[r]^{pr_{1}^{*}}\ar@{.>}[d]& \cV ect^{G}\( M\times E^{f}_{G}\)\\
\cV ect(M'=M/G)\ar[r]_{{pr'}^{*}_{1}}^{\approx^{\otimes\cQ}} & \cV ect\(\(M\times E^{f}_{G}\)/G\)\ar[u]_{q^{*}}^{\approx}
}
}$$

In the case of classical signature we can apply the following theorem:

\begin{theorem}
The quotient space ${X=M/G}$ satisfies the Poincare duality
in homology with rational coefficients:
$${D:H^{k}(X;\Q)\mapr{\approx}H_{n-k}(X;\Q)}$$
\end{theorem}

It is easy to notice that the Hirzebruch genus of the manifold $M$
can be represented by an invariant differential form with respect
to the action of the group $G$. Therefore, one can define the Hirzebruch genus
of the quotient space $M/G$.

\begin{conjecture}
For the signature ${\sign M/G=\sign H^{*}(M/G;\cQ)\in \mathbb{Z}}$
one has an analogue of the Hirzebruch formula:
$${
\sign M/G=\<L(M), [M/G]\>\in \cQ.
}$$
\end{conjecture}

\textbf{For a $C^{*}$--algebra $\cA$:}
In the case of a $C^{*}$--algebra $\cA=C^{*}[G]$, a similar
commutative diagram holds:
$${
\xymatrix{
{\cV ect}^{G}(M; \cA)\ar[r]^-{pr_{1}^{*}}\ar@{.>}[d]& \cV ect^{G}\( M\times E^{f}_{G}; \cA\)\\
\cV ect(M/G; \cA)\ar[r]_-{{pr'}^{*}_{1}}^-{\approx^{\otimes\cQ}} & \cV ect\(\(M\times E^{f}_{G}\)/G; \cA\)\ar[u]_{q^{*}}^{\approx}
}
}$$

%\textbf{For a $C^{*}$--algebra $\cA$:}
Define ${\xi_{\cA}\in {\cV ect}^{G}(M; \cA)}$ by ${
\xi_{\cA}\bydef M\times\cA.}$ Then, there is a number ${n\in \Z}$
and a bundle ${\eta_{\cA}\in\cV ect(M/G; \cA) }$ such that
${
q^{*}{pr'}_{1}^{*}(\eta_{\cA})=n\cdot pr_{1}^{*}(\xi_{\cA})
}$.

\begin{conjecture}{Conjecture for the $C^{*}$--algebra $\cA$:}
For the noncommutative signature ${\sign M\in K(\cA)}$ one has
$${
\sign M=\<\frac{1}{n}L(M)\ch_{\cA}(\eta_{\cA}), [M/G]\>\in K(\cA)\otimes\Q.
}$$
\end{conjecture}

\section{The Connor--Floyd's fix-point construction}

To clarify the bordism concept for proper action one can apply so
called the Conner-Floyd construction for fixed points \cite{Conner-1964}.
This construction reduces the problem of describing these bordisms to two simpler
problems: a) description of the fixed-point sets (or, more generally,
the stationary point sets), which happens to be a submanifold attached
with the structure of its normal bundle and the action of the same group $G$,
however, this action could have only stationary points of lower rank; b) description of the bordisms of lower rank with an action of the group $G$.

Let
$${\cH(M)\bydef \{H_{x}: x\in M, \quad H_{x}\neq \{1\}\}.}$$
be the category of isotropic groups.

Let $[H]$ be the family of all subgroups
conjugated to the subgroup
${H\subset G}$. Let ${\rF}$ be a family of subgroups
that is invariant with respect to conjugation and closed
under inclusions. Denote by
${\max\{\rF\}}$  the subfamily  of maximal objects in  ${\rF}$.

The family $\rH_{0}$ is already closed under conjugation and
the subgroups of its elements can be added to it to obtain the
whole family of subgroups of $G$ having nonempty fixed points set.
So, put ${\rH_{0}\bydef \max\{\cH(M)\}}$, ${\rH'_{0}\bydef \cH(M)\setminus\rH_{0}}$,
and
${\rH_{k}\bydef \max\{\rH'_{k-1}\}}$,
${\rH'_{k}\bydef \rH'_{k-1}\setminus\rH'_{k-1}}$.
All the families
$$
\rH_{0} \succ \rH_{1}\succ \rH_{2} \succ \dots \succ \rH_{k} \succ \dots
$$
are invariant with respect to conjugation.

The set of stationary points  ${M^{\rH_{0}}}$ is a submanifold invariant
with respect to the action of the group ${G}$.
A tubular invariant neighborhood
${U_{0}\supset M^{\rH_{0}}}$\linebreak is diffeomorphic to a normal vector bundle
${\xi_{0}}$ with base ${M^{\rH_{0}}}$,\linebreak
${\dim \xi_{0} + \dim M^{\rH_{0}} = \dim M = n}$.

So one has an exact sequence of equivariant bordisms
$${
\begin{array}{ccccccc}
\cdots\mapr{}&{}_{n}\Omega_{G}(\rH'_{0})&
  \mapr{inclusion}&{}_{n}\Omega_{G}(\rH_{0}\sqcup\rH'_{0})&\mapr{fixed point}\\\\
  \mapr{fixed point}&{}_{n}\Omega_{G}^{v}(\rH_{0},\rH'_{0})
  & \mapr{\partial} & {}_{n-1}\Omega_{G}(\rH'_{0}) & \mapr{}
  \cdots
\end{array}
}$$

Similar for any integer ${k}$ one has also a exact sequence:
$${
\begin{array}{ccccccc}
\cdots\mapr{}&{}_{n}\Omega_{G}(\rH'_{k})&
  \mapr{inclusion}&{}_{n}\Omega_{G}(\rH'_{k-1})&\mapr{fixed point}\\\\
  \mapr{fixed point}&{}_{n}\Omega_{G}^{v}(\rH_{k},\rH'_{k})
  & \mapr{\partial} & {}_{n-1}\Omega_{G}(\rH'_{k}) & \mapr{}
  \cdots
\end{array}
}$$

Here, the group $ {}_{n}\Omega_{G}^{v}(\rH_{k},\rH'_{k})$ is
generated by manifolds equipped with the structure of their
normal bundle with action of the group $G$

For an element
$H\in \rH_{0}$ the restriction ${\xi_{H}}$ of the bundle ${\xi_{0}}$
to the fixed points set $M^H$ has a fiberwise action
$N_G(H)\times {\xi_{H}}\rightarrow {\xi_{H}}$
of the normalizer $N_G(H)$ such that the diagram
\ff{
%\ff{
\begin{array}{ccc}
  N_G(H)  \times{\xi_{H}} &\mapr{} &{\xi_{H}}\\
  \mapd{}           &        &\mapd{}\\
  N_G(H)/H\times M^H  &\mapr{} &M^H \\
\end{array}
}is commutative, and the lower action is free\footnote{
This is called quasi-free action of the group $N_G(H)$ in
\cite{Levine-2008}.}.

For an element
$gHg^{-1}\in \rH_{0}$ in the conjugation class of $H$
the bundles ${\xi_{H}}$ and ${\xi_{gHg^{-1}}}$ are
equivariantly isomorphic.

And, by the maximality of the elements
$H,K\in \rH_{0}$ the corresponding fixed points sets
$M^H$ and $M^K$ and are disjoint.

Therefore, the group $ {}_{n}\Omega_{G}^{v}(\rH_{k},\rH'_{k})$
splits as finite direct sum
$$ {}_{n}\Omega_{G}^{v}(\rH_{k},\rH'_{k})\approx
 \bigoplus_{[H]\in\rH_{k}/G}{}_{n}\Omega_{G}^{v}([H]).$$

So, the calculation of this groups reduce to the
description of the classifying space for equivariant vector bundles
for the case of quasi-free action of the group
$G$ on the base.

\section{Classifying space for vector $G$-bundles with quasi-free proper
action of a discrete group}

This problem naturally arises from the Connor-Floyd's
\cite{Conner-1964} description
of bordisms with the action of a group
$G$ using the so-called fix-point construction
% при помощи так называемой фикспоинт-конструкции.

Consider a finite subgroup  $H<G$ and the
set $M^{H}$ of fixed points with respect to the action of the subgroup $H$.
If the subgroup $H\in\max\{\cH(M)\}$, then its conjugate groups $H'=gHg^{-1}$
either lead to the same fixed point set, i.e. $M^{H}=M^{H'}$, or these
sets do not intersect, i.e. $M^{H}\cap M^{H'}=\emptyset$. In any case, the union
$$M'=\bigcup\limits_{g\in G}M^{gHg^{-1}}$$ is a submanifold, whose regular tubular
neighborhood is equivariantly diffeomorphic to a vector $G$-bundle $\xi$, and the
action of the group $G$ over the total space of this bundle does not have $H$ or any
of its conjugates as stationary subgroups. So, the submanifold $M'$ decomposes into its components
$$
M'=\bigsqcup\limits_{[g]\in G/N(H)}g(M^{H}),
$$
where $N(H)$ is the normalizer of the subgroup $H$ and, over the set $M^{H}$, acts the factor group $N(H)/H$.

This shows that, in order to describe the structure of the
fixed points sets using the Connor-Floyd's method, it is enough to find a
description of the equivariant vector bundles with a particular kind of action
of the group $G$ (or its normalizer $N(H)$) over the base of the bundle:
the so-called quasi-free action. Previously, by different authors (like, for example, \cite{Atiyah-1967}), there were considered only free and trivial actions of a
 discrete group $G$. The case of an arbitrary action was also studied in the work
(\cite[7.2]{Atiyah-Segal-1965}), although not for bundles, but for the
$K$-functor generated by them, and only for compact groups.
For these reasons, we believe that the problem of describing
the equivariant vector bundles for proper quasi-free action of discrete
groups shall be of great interest.

A brief exposition and several preliminary results were published in 
\linebreak \cite{MishQuit},\cite{Quit},\cite{MishQuit-1},\cite{[2]},\cite{[3]},\cite{[1]},\cite{[4]}.
%The first author was partially supported by the grants RFFI No.08-01-00034-a, No. 10-01-92601-KO-a and for the project of the Ministry of Education of
%the Russian Federation No. 2.1.1/5031

\subsection{The setting of the problem}

Lets  $\xi$ be an $G$-equivariant vector bundle with base $M$.
\beq{1}
\begin{array}{c}
          \xi \\
          \mapd{} \\
          M \\
        \end{array}
\eeq
Lets $H\lhd G$ be a normal finite subgroup. Assume that the action of the group
$G$  over the base $M$ reduces to the factor group $G_{0}=G/H$:
\beq{2}
\xymatrix{
  G\times M \ar[dr]\ar[dd]&\\
  &M .\\
  G_{0}\times M \ar[ur]& \\
}\eeq
suppose, additionally,
that the action $G_{0}\times M\mapr{} M$ is free and there is no more fixed points
of the action of the group $H$ in the total space of the bundle $\xi$.

So, we have the following commutative diagram

\feq{3}{
\begin{array}{ccc}
  G\times\xi &\mapr{}&\xi\\
  \mapd{} &&\mapd{}\\
  G_{0}\times M&\mapr{}&M \\
\end{array}
}

\begin{definition}
As in \cite[p. 210]{Levine-2008}, we shall say that the described
action of the group $G$ is \textit{quasi-free}
over the base with normal \textit{stationary} subgroup $H$.
\end{definition}

Reducing the action to the subgroup $H$, we obtain the simpler diagram:
\feq{4}{
\begin{array}{ccc}
  H\times\xi &\mapr{}&\xi\\
  \mapd{} &&\mapd{}\\
   M&=&M \\
\end{array}
}

Let $\rho_{k}: H\mapr{}\U(V_{k})$ be the series of all the irreducible
(unitary) representation of the finite group $H$. Then the $H$-bundle $\xi$
can be presented as the finite direct sum:
\feq{5}{
\xi\approx \bigoplus_{k}\(\xi_{k}\bigotimes V_{k}\),
}
where the action of the group  $H$ over the bundles $\xi_{k}$ is trivial,  $V_{k}$
denotes the trivial bundle with fiber $V_{k}$ and with  fiberwise
action of the group $H$, defined using the linear representation $\rho_{k}$.

\begin{lemma} The group $G$ acts on every term of the sum (\ref{5})
separately, i.e. every summand $\xi_{k}\bigotimes V_{k}$ is invariant
under the action of the group $G$.
\end{lemma}

\proof\, Consider now the action of the group $G$ over the total space
of the bundle $\xi$. Fix a point $x\in M$. The action of the element
$g\in G$ is fiberwise, and maps the fiber $\xi_{x}$ to the fiber
$\xi_{gx}$: \ff{ \Phi(x,g):\xi_{x}\mapr{}\xi_{gx}. }

 Also, for a par of elements  $g_{1}, g_{2}\in G$ we have:

\feq{5.1}{
\Phi(x, g_{1}g_{2})=\Phi\(g_{2}x, g_{1}\)\circ\Phi\(x, g_{2}\),
}

\ff{
\begin{array}{ccccc}
 \Phi(x, g_{1}g_{2}):     \xi_{x}& \mapr{\Phi\(x, g_{2}\)}&\xi_{g_{2}x}&
      \mapr{\Phi\(g_{2}x, g_{1}\)} & \xi_{g_{1}g_{2}x}\\
\end{array}
}

In particular, if $g_{2}=h\in H\lhd G$, then $g_{2}x=hx=x$. So,
\ff{
\begin{array}{ccccc}
 \Phi(x, gh):     \xi_{x}& \mapr{\Phi\(x, h\)}&\xi_{x}&
      \mapr{\Phi\(x, g\)} & \xi_{gx}\\
\end{array}
} Analogously, if $g_{1}=h\in H<G$, then $g_{1}gx=hgx=gx$. So \ff{
\begin{array}{ccccc}
 \Phi(x, hg):     \xi_{x}&
      \mapr{\Phi\(x, g\)} & \xi_{gx}& \mapr{\Phi\(gx, h\)}\xi_{gx}\\
\end{array}
}
The operator $\Phi\(x, h\)$ does not depends on the point $x\in M$,

\ff{
\Phi(x, h)=\Psi(h):\bigoplus_{k}\(\xi_{k,x}\bigotimes V_{k}\)\mapr{}
\bigoplus_{k}\(\xi_{k,x}\bigotimes V_{k}\),
}here, since the action of the group  $H$ is given over every space
$V_{k}$ using pairwise different irreducible representations $\rho_{k}$,
we have
\ff{
\Psi(h)=\bigoplus_{k}\(\id\bigotimes \rho_{k}(h)\).
}

In this way, we obtain the following relation:
\feq{6}{
\Phi(x, gh)=\Phi(x, g)\circ\Psi(h)=
\Phi(x, ghg^{-1}g)=\Psi(ghg^{-1})\circ\Phi(x, g).
}

Lets write the operator  $\Phi(x, g)$ using matrices to decompose the space
 $\xi_{x}$ as the direct sum
\ff{
\xi_{x}=\bigoplus_{k}\(\xi_{k,x}\bigotimes V_{k}\):
}
\feq{6.1}{
\Phi(x, g)=\(
\begin{array}{cccc}
   \Phi(x, g)_{1,1}& \cdots & \Phi(x, g)_{k,1} & \cdots \\
  \vdots & \ddots & \vdots &  \\
  \Phi(x, g)_{1,k} & \cdots & \Phi(x, g)_{k,k} & \cdots\\
  \vdots &  & \vdots & \ddots \\
\end{array}
\)
}

If $k\neq l$ then $\Phi(x, g)_{k,l}=0$, i.e. the matrix $\Phi(x, g)$ its diagonal,
\ff{
\Phi(x, g)=\bigoplus_{k}\Phi(x, g)_{k,k}:
\bigoplus_{k}\(\xi_{k,x}\bigotimes V_{k}\)\mapr{}
\bigoplus_{k}\(\xi_{k,gx}\bigotimes V_{k}\),
}
\ff{
\Phi(x, g)_{k,k}:
\(\xi_{k,x}\bigotimes V_{k}\)\mapr{}
\(\xi_{k,gx}\bigotimes V_{k}\),
}
as it was required to prove. $\hfill\blacksquare$

\subsection{Description of the particular case $\xi=\xi_{0}\bigotimes V$}
Here we will consider the particular case of a $G$-vector bundle $\xi=\xi_{0}\otimes V$ with base $M$.
\ff{
\begin{array}{c}
          \xi_{0}\otimes V \\
          \mapd{} \\
          M \\
        \end{array}
}
where  the action of the group $G$ is quasi-free over the base with finite normal stationary subgroup  $H\lhd G$.

We will assume that the group $H$ acts trivially over the bundle $\xi_0$.
By $V$ we denote the trivial bundle with fiber $V$ and with fiberwise action
of the group $H$ given by an irreducible linear representation $\rho$.

\begin{definition} A \textit{canonical model} for the fiber in a $G$-bundle
$\xi=\xi_0\bigotimes V$ with fiber $F\otimes V$ is the bundle over
one orbit homeomorphic to $G_0$ and with fiber $(F\otimes V)$, i.e.
the cartesian product $\cW=G_0\times \(F\otimes V\)$ with the natural projection
$$\cW=G_0\times \(F\otimes V\)\mapr{} G_{0}$$
and fiberwise action of the group $G$ \ff{\begin{array}{ccc}
  G\times \cW &\mapr{\phi}& \cW\\
  \mapd{} &&\mapd{}\\
  G\times G_{0}&\mapr{\mu}&G_{0} \\
\end{array}
}i.e.
\ff{
\begin{array}{ccc}
  G\times \(G_0\times \(F\otimes V\)\) &\mapr{\phi}& G_0\times \(F\otimes V\)\\
  \mapd{} &&\mapd{}\\
  G\times G_{0}&\mapr{\mu}&G_{0} \\
\end{array}
}where
$\mu$ denotes the natural left action of $G$ on its quotient $G_0$, and
\ff{
\phi([g],g_{1}):[g]\times \(F\otimes V\) \to [g_1g]\times \(F\otimes V\)
}is
given by the formula
\feq{15.1a}{
\begin{array}{ll}
\phi([g],g_{1})
=
\id\otimes
\rho(u(g_{1}g)u^{-1}(g)).
    \end{array}
}where
\ff{
u:G\mapr{}H
}
is a homomorphism of right $H$-modules by multiplication,
i.e.
\ff{
u(gh)=u(g)h,\quad  u(1)=1, \quad g\in G, h\in H.
}
\end{definition}

\begin{lemma} The definition (\ref{15.1a}) of the action of  \f{G} is
well-defined.
\end{lemma}

\textsc{Proof.} It is enough to prove that
that  a) the formula  (\ref{15.1a}) defines an associative action, i.e.
\ff{
\begin{array}{ll}
\phi([g],g_2g_{1})
=\phi([g_1g],g_2)\circ\phi([g],g_{1}),
    \end{array}
}for every $g\in G$, $g_{1}\in G$, $g_{2}\in G$ and

b) that the formula (\ref{15.1a}) does not depends on the chosen representative $gh\in [g]$:
\ff{\id\otimes \rho(u(g_{1}g)u^{-1}(g))=\id\otimes \rho(u(g_{1}gh)u^{-1}(gh))}for
 every $g\in G$, $g_{1}\in G$ and $h\in H$.

In fact,
\ff{
\begin{array}{ll}
\phi([g],g_2g_{1})
=
\id\otimes \rho(u(g_2g_{1}g)u^{-1}(g))=\\\\
\id\otimes \rho(u(g_2g_{1}g)u(g_{1}g)u^{-1}(g_{1}g)u^{-1}(g))=\\
=\id\otimes \rho(u(g_2g_{1}g)u(g_{1}g))\circ\id\otimes \rho(u^{-1}(g_{1}g)u^{-1}(g))=\\\\
=\phi([g_1g],g_2)\circ\phi([g],g_{1}),
    \end{array}
}what proves a), and,  recalling the equation $u(gh)=u(g)h$ for every $g\in G$ and $h\in H$, it is clear that
\ff{u(g_{1}gh)u^{-1}(gh)=u(g_{1}g)hh^{-1}u^{-1}(g)=u(g_{1}g)u^{-1}(g),}
 which implies b). $\hfill\blacksquare$

In general, a canonical model depends on the choice of the homomorphism $u:G\mapr{}H$,
fortunately, any two such canonical models turn out to be isomorphic, as it
is shown in the following lemma.

\begin{lemma}
For different choices of homomorphisms $u:G\mapr{}H$ and $u':G\mapr{}H$ the corresponding canonical models $\cW$ and $\cW'$ are equivariantly isomorphic,
i.e. there is a fiberwise equivariant homomorphism
$$
\begin{array}{ccc}
\cW&\mapr{\psi}&\cW' \\
\mapd{}&&\mapd{}\\
G_{0}&=&G_{0}
\end{array}
$$
\end{lemma}
\proof\,
It is required to construct an equivariant
isomorphism

\ff{
\begin{array}{ccc}
  G_0\times \(F\otimes V\) &\mapr{\psi}& G_0\times \(F\otimes V\)\\
  \mapd{} &&\mapd{}\\
 G_{0}&=&G_{0} \\
\end{array}
}
We define  $\psi$ by the formula
$$
\psi([g]):[g]\times (F\otimes V)\mapr{}[g]\times (F\otimes V)
$$

\begin{equation}\label{10}
\psi([g])=\id\otimes\rho\(u'(g)u^{-1}(g)\)
\end{equation}

The formula (\ref{10}) is well define, since a) does not depends on the choice of
the representative in the conjugation class
$g\in[g]$, also, it defines an equivariant mapping between canonical models,
i.e. the diagram

\ff{
\begin{array}{ccc}
  [g]\times \(F\otimes V\) &\mapr{\phi(g_{1},[g])}& [g_{1}g]\times \(F\otimes V\)\\
  \mapd{\psi([g])} &&\mapd{\psi([g_{1}g])}\\
  {[g]\times \(F\otimes V\)} &\mapr{\phi'(g_{1},[g])}& {[g_{1}g]\times \(F\otimes V\)}
\end{array}
}
commutes.
$\hfill\blacksquare$

Lets consider, over the base $M$, an atlas of
equivariant charts $\{O_{\alpha}\}$,
\ff{
M=\bigcup_{\alpha}O_{\alpha},
}
\ff{
[g]O_{\alpha}=O_{\alpha}, \qquad [g]\in G_{0}.
}

\begin{lemma}\label{l4}
There exists an atlas  $\{O_{\alpha}\}$ fine enough,
such that every chart $O_{\alpha}$ can be presented as a disjoint union
of its subcharts:
\ff{
O_{\alpha}= \bigsqcup_{[g]\in G_{0}}[g]U_{\alpha},
}
where the translates of the chart $U_{\alpha}$ are pair-wise disjoint,
i.e.
\ff{[g]U_{\alpha}\cap [g']U_{\alpha}=\emptyset \hbox{ for }  [g]\neq [g'],} and,
for different indices  $\alpha\neq\beta$,  a nonempty intersection
$U_{\alpha}\cap [g_{\alpha\beta}]U_{\beta}\neq \emptyset$,
can occur only for an unique element  $[g_{\alpha\beta}]\in G_{0}$.
\end{lemma}

In other words,
if $[g]\neq[g_{\alpha\beta}]$,
then $U_{\alpha}\cap[g]U_{\beta}=\emptyset$.

This means, that every chart  $O_{\alpha}$ is homeomorphic
to the cartesian product
$$
O_{\alpha}\approx U_{\alpha}\times G_{0},
$$
and also the intersection of two charts $O_{\alpha}\cap O_{\beta}$
can be presented as the cartesian product
$$
O_{\alpha}\cap O_{\beta}\approx (U_{\alpha}\cap [g_{\alpha\beta}]U_{\beta})\times G_{0}.
$$

In terms of the introduced notations, the following theorem constitutes
an analog of the proposition about locally trivial bundles (conf., for example, \cite{Mishchenko-1998}).

\begin{theorem}
The bundle  $\xi= \xi_{0}\bigotimes V$ is locally homeomorphic
to the cartesian product of some chart $U_\alpha$ by the canonical
model.  More precisely, for a fine enough atlas, there exist
$G$-equivariant trivializations
\feq{14.a}{\psi_\alpha:O_\alpha\times \(F\otimes
V\)\to\xi|_{O_\alpha}}where \ff{O_\alpha \times \(F\otimes
V\)\approx U_\alpha \times \(G_0\times \(F\otimes V\)\)}and there is
a commutative diagram 
\feq{trivial}{
\begin{array}{ccc}\xi|_{O_\alpha}&\mapr{g}&\xi|_{O_\alpha}\\
                               \mapu{\psi_\alpha} &&\mapu{\psi_\alpha}\\
                               U_\alpha\times \(G_0\times \(F\otimes V\)\)
                               &\mapr{\id\times \phi(g)}&U_\alpha\times \(G_0\times \(F\otimes V\)\)\\
\end{array}
}where $g\in G$, \f{\id:U_\alpha\to U_\alpha,} and $\phi(g)$ denotes
the canonical action.
\end{theorem}

\proof\, Consider an atlas fine enough as to have the properties stated
in the lemma \ref{l4}. We will construct the trivialization (\ref{14.a})
from an arbitrary $H$ equivariant trivialization \ff{\psi_\alpha:U_\alpha\times
\(F\otimes V\)\to\xi|_{U_\alpha}}in such a way, that the diagram
\begin{equation}\label{13}
\begin{array}{ccc}\xi|_{U_\alpha}&\mapr{g}&\xi|_{[g]U_\alpha}\\
\mapu{\psi_{\alpha,1}} &&\mapu{\psi_{\alpha,[g]}}\\
U_\alpha\times \(F\otimes V\)
&\mapr{\id\times\phi(g,[1])}&[g]U_\alpha\times  \(F\otimes V\)\\
\end{array}
\end{equation}
commutes for every $g\in [g]$, where the left, upper and lower arrows
are given and the right one is defined by the commutativity of the diagram.
It is enough to verify that the map $\psi_{\alpha,[g]}$ is well-defined,
i.e. that it does not depend on the conjugation class $[g]$,
$$
\psi_{\alpha,[g]}=\psi_{\alpha,[hg]}, \quad h\in H.
$$

In fact, since the isomorphism $\psi_{\alpha,1}$ is  $H$-equivariant,
so the following diagram commutes:
\begin{equation}\label{14}
\begin{array}{ccc}\xi|_{U_\alpha}&\mapr{h}&\xi|_{[g]U_\alpha}\\
\mapu{\psi_{\alpha,1}} &&\mapu{\psi_{\alpha,1}}\\
U_\alpha\times \(F\otimes V\)
&\mapr{\id\times\phi(h,[1])}&U_\alpha\times  \(F\otimes V\)\\
\end{array}
\end{equation}
Lets glue the diagrams  (\ref{13},\ref{14}) together

\begin{equation}\label{15}
\begin{array}{ccccc}\xi|_{U_\alpha}&\mapr{h}&\xi|_{[g]U_\alpha}&\mapr{g}&\xi|_{[g]U_\alpha}\\
\mapu{\psi_{\alpha,1}} &&\mapu{\psi_{\alpha,1}}&&\mapu{\psi_{\alpha,[g]}}\\
U_\alpha\times \(F\otimes V\)
&\mapr{\id\times\phi(h,[1])}&U_\alpha\times  \(F\otimes V\)&\mapr{\id\times\phi(g,[1])}&[g]U_\alpha\times  \(F\otimes V\)\\
\end{array}
\end{equation}
we obtain the following diagram:
\begin{equation}\label{16}
\begin{array}{ccc}\xi|_{U_\alpha}&\mapr{gh}&\xi|_{[g]U_\alpha}\\
\mapu{\psi_{\alpha,1}} &&\mapu{\psi_{\alpha,[g]}}\\
U_\alpha\times \(F\otimes V\)
&\mapr{\id\times\phi(gh,[1])}&[g]U_\alpha\times  \(F\otimes V\)\\
\end{array}
\end{equation}
since
$$
\phi(g,[1])\cdot\phi(h,[1])=\phi(gh,[1]),
$$
but this means that the mapping
$\psi_{\alpha,[g]}$ is well-defined, as it was
required to prove.

 $\hfill\blacksquare$

By \f{\mathrm{Aut}_G\(G_0\times \(F\otimes V\)\)} we denote the group of equivariant automorphisms
of the space $G_0\times \(F\otimes V\)$ as a vector $G$-bundle with base $G_0$, fiber
$F\otimes V $ and canonical action of the group $G$.

\begin{corollary}\label{transitionfunctions} The transition functions on the intersection
\f{ \(U_\alpha\times G_0\) \cap \(U_\beta\times G_0\)=
\(U_\alpha\cap [g_{\alpha\beta}]U_\beta\)\times G_0},
i.e. the homomorphisms \f{\Psi_{\alpha\beta}=\psi_{\beta}^{-1}\psi_{\alpha}}
on the diagram
\feq{15.2a}{
\begin{array}{ccc}
  \(U_\alpha\cap  [g_{\alpha\beta}]U_\beta\)\times \(G_0\times \(F\otimes V\)\)
  &\mapr{\Psi_{\alpha\beta}}& \(U_\alpha\cap  [g_{\alpha\beta}]U_\beta\)\times \(G_0\times \(F\otimes V\)\)\\
  \mapd{} &&\mapd{}\\
   \(U_\alpha\cap  [g_{\alpha\beta}]U_\beta\)\times G_0
  &\mapr{\id}& \(U_\alpha\cap [g_{\alpha\beta}]U_\beta\)\times G_0\\
\end{array}
}are equivariant with respect to the canonical action of the group
$G$ over the product of the base by the canonical model, i.e.
\ff{\Psi_{\alpha\beta}(x)\circ
\phi(g_1,[g])=\phi(g_1,[g])\circ\Psi_{\alpha\beta}(x)}for every
$x\in U_\alpha\cap  [g_{\alpha\beta}]U_\beta,\;g_1\in G,\; [g]\in
G_0$. In other words, \ff{\Psi_{\alpha\beta}(x)\in
\mathrm{Aut}_G\(G_0\times \(F\otimes V\)\).}
\end{corollary}

Now we give a more accurate description of the group
\f{\mathrm{Aut}_G\(G_0\times \(F\otimes V\)\)}.
By definition, an element of the group \f{\mathrm{Aut}_G\(G_0\times \(F\otimes V\)\)}
is an equivariant mapping $\mathbf{A}^{a}$, such that the pair $(\mathbf{A}^{a},a)$
defines a commutative diagram
\ff{
\begin{array}{ccc}
   \(G_0\times \(F\otimes V\)\) &\mapr{\mathbf{A}^{a}}& G_0\times \(F\otimes V\)\\
  \mapd{} &&\mapd{}\\
   G_{0}&\mapr{a}&G_{0}, \\
\end{array}
}
which commutes with the canonical action, i.e.
the map \f{a\in \mathrm{Aut}_G(G_0)} satisfies the condition
\ff{a\in \mathrm{Aut}_G(G_0)\approx G_0,\quad a[g]=[ga],\;[g]\in G_0,}
and the mapping   $\mathbf{A}^{a}=(A^{a}[g])_{[g]\in G_0}$,

\ff{
A^{a}[g]: [g]\times(F\otimes V)\to [ga]\times (F\otimes V)
}
satisfies a commutation condition with respect to the action of the group \f{G}:
\ff{
\begin{array}{ccc}
  [g]\times(F\otimes V) & \mapr{A^{a}[g]} & [ga]\times(F\otimes V) \\
  \mapd{\phi(g_{1},[g])} &  & \mapd{\phi(g_{1},[ga])} \\
  {[g_{1}g]\times(F\otimes V)} & \mapr{A^{a}[g_{1}g]} & [g_{1}ga]\times(F\otimes V)
\end{array} \quad ,
}

\feq{16.e}{\phi(g_1,[ga])\circ A^{a}[g]=A^{a}[g_{1}g]\circ\phi(g_1,[g])}
 i.e.
 \feq{16.b}{
 (\id\otimes \rho(u(g_{1}ga)u^{-1}(ga)))A^{a}[g]=
 A^{a}[g_1g](\id\otimes\rho(u(g_{1}g)u^{-1}(g)))
 }
where  $ [g]\in G_0, \quad g_1\in G$.

\begin{lemma}\label{shortexact} It takes place an exact sequence of groups
\feq{suc}{1\to GL(F)
\mapr{}\mathrm{Aut}_G\(G_0\times \(F\otimes V\)\)
\mapr{} G_0\to 1}
\end{lemma}

\proof\, To define a projection
\ff{pr:\mathrm{Aut}_G\(G_0\times \(F\otimes V\)\) \mapr{} G_0}we
send the fiberwise map
\ff{\mathbf{A}^{a}:G_0\times \(F\otimes V\) \mapr{} G_0\times \(F\otimes V\)}to
its restriction over the base
\f{a:G_0\to G_0}, i.e. $a\in \mathrm{Aut}_G(G_0)\approx G_0$. So, this is a well-defined homomorphism.

We need to show that $pr$ is an epimorphism and that its kernel is isomorphic to $GL(F)$.
Lets calculate the kernel, i.e. the set of mappings of the form $A^{1}$.
The equation (\ref{16.b}) in this case gives
 \feq{17.a}{
 (\id\otimes \rho(u(g_{1}g)u^{-1}(g)))A^{1}[g]=
 A^{1}[g_1g](\id\otimes\rho(u(g_{1}g)u^{-1}(g)))
 }

In the case  $g_{1}=h\in H$, we obtain

 \ff{
 (\id\otimes \rho(u(hg)u^{-1}(g)))A^{1}[g]=
 A^{1}[g](\id\otimes\rho(u(hg)u^{-1}(g)))
 }
Since the representation \f{\rho} is irreducible,
by Schur's lemma, we have
\ff{
A^{1}[g] = B^{1}[g]\otimes \id.}

On the other side, assuming in (\ref{17.a}) that \f{g=1}, we have

 \ff{
 (\id\otimes \rho(u(g)))A^{1}[1]=
 A^{1}[g](\id\otimes\rho(u(g))),
 }
 i.e.

  \ff{
 (\id\otimes \rho(u(g)))(B^{1}[1]\otimes \id)=
 (B^{1}[g]\otimes \id)(\id\otimes\rho(u(g))),
 }
 or
   \ff{
 (B^{1}[g]\otimes \id)=
 (B^{1}[1]\otimes \id).
 }

 So, the kernel \f{\ker pr} is isomorphic to the group \f{GL(F)}.

 In the generic case, i.e. \f{[a]\neq 1},
 we can compute the operator \f{A^{a}[g]} in terms of its value at the identity \f{A^{a}[1]}
 from the formula (\ref{16.b}):
assuming  $g=1$, we obtain (changing $g_{1}$ by $g$):

  \feq{16.c}{
 (\id\otimes \rho(u(ga)u^{-1}(a)))A^{a}[1]=
 A^{a}[g](\id\otimes\rho(u(g))),
 }

 i.e.

  \feq{16.d}{
A^{a}[g]=
(\id\otimes \rho(u(ga)u^{-1}(a)))A^{a}[1](\id\otimes\rho(u^{-1}(g))),
 }

 Therefore, the operator is completely defined by its value
\ff{
A^{a}[1]: [1]\times(F\otimes V)\to [a]\times (F\otimes V)
}
at the identity $g=1$.

Now we describe the operator $A^{a}[1]$ in terms of the representation
$\rho$ and its properties.

We have a commutation rule with respect to the action of the subgroup \f{H}:
\ff{
\begin{array}{ccc}
  [1]\times(F\otimes V) & \mapr{A^{a}[1]} & [a]\times(F\otimes V) \\
  \mapd{\phi(h,[1])} &  & \mapd{\phi(h,[a])} \\
  {[1]\times(F\otimes V)} & \mapr{A^{a}[1]} & [a]\times(F\otimes V)
\end{array} \quad ,
}
Equivalently

\ff{
A^{a}[1]\circ\phi(h,[1])=\phi(h,[a])\circ A^{a}[1],
}

i.e.

\ff{
A^{a}[1]\circ(\id\otimes \rho(h))=(\id\otimes\rho(g'^{-1}(a)hg'(a)))\circ A^{a}[1],
}

i.e.

\ff{
A^{a}[1]\circ(\id\otimes \rho(h))=(\id\otimes\rho_{g'(a)}(h))\circ A^{a}[1].
}

The last equation means that the operator $A^{a}[1]$ should  permute these
representations, or equivalently, such an operator exists only when the
representations $\rho$ and $\rho_{g'(a)}$ are equivalent. Recalling
the commutation rule (\ref{6}), we see that this is the case we are been
considering.

Thus, if the representations  $\rho$ and $\rho_{g}$ are equivalent,
 we have an (invertible) splitting operator \f{C(g)},
satisfying the equation
\feq{29-a}{
\rho_{g}(h)=\rho\(g^{-1}hg\)=C(g)\rho(h)C^{-1}(g).
}for every $g\in G$.
The operator  \f{C(g)} is defined up to multiplication by a scalar operator
\f{\mu_{g}\in\SS^{1}\subset\C^{1}}.

So

\ff{
A^{a}[1]\circ(\id\otimes \rho(h))=(\id\otimes C(g'(a))\circ\rho(h)\circ C^{-1}(g'(a)))\circ A^{a}[1],
}
or

\ff{
(\id\otimes C^{-1}(g'(a)))\circ A^{a}[1]\circ(\id\otimes \rho(h))=(\id\otimes\rho(h))\circ (\id\otimes C^{-1}(g'(a)))\circ A^{a}[1],
}

Then, by the Schur's lemma,

\ff{
(\id\otimes C^{-1}(g'(a)))\circ A^{a}[1]=B^{a}[1]\otimes\id,
}
i.e.
\ff{
A^{a}[1]=B^{a}[1]\otimes C(g'(a)),
}

Using the formula (\ref{16.d}), we obtain
  \ff{
A^{a}[g]=
(\id\otimes \rho(u(ga)u^{-1}(a)))(B^{a}[1]\otimes C(g'(a)))(\id\otimes\rho(u^{-1}(g))),
 }

i.e.

  \feq{16.f}{
A^{a}[g]=
B^{a}[1]\otimes (\rho(u(ga)u^{-1}(a))\circ C(g'(a))\circ\rho(u^{-1}(g))).
 }

This means, that by defining the matrix $B^{a}[1]$, it is possible to obtain all the operators
$A^{a}[g]$ satisfying the equation (\ref{16.d}).

It remains to verify the commutation rule (\ref{16.b}), i.e. in the formula

 \ff{
 (\id\otimes \rho(u(g_{1}ga)u^{-1}(ga)))A^{a}[g]=
 A^{a}[g_1g](\id\otimes\rho(u(g_{1}g)u^{-1}(g)))
 }
we substitute the expression (\ref{16.f}):
\ff{
\begin{array}{c}
  (\id\otimes \rho(u(g_{1}ga)u^{-1}(ga)))\circ(
  B^{a}[1]\otimes (\rho(u(ga)u^{-1}(a))\circ C(g'(a))\circ\rho(u^{-1}(g)))) =\\\\
  =(B^{a}[1]\otimes (\rho(u(g_{1}ga)u^{-1}(a))\circ C(g'(a))\circ\rho(u^{-1}(g_{1}g))))
  \circ(\id\otimes\rho(u(g_{1}g)u^{-1}(g)))
\end{array}
}

that is
\ff{
\begin{array}{c}
  B^{a}[1]\otimes \rho(u(g_{1}ga)u^{-1}(ga)))\circ
  (\rho(u(ga)u^{-1}(a))\circ C(g'(a))\circ\rho(u^{-1}(g)))) =\\\\
  =B^{a}[1]\otimes (\rho(u(g_{1}ga)u^{-1}(a))\circ C(g'(a))\circ\rho(u^{-1}(g_{1}g))))
  \circ(\rho(u(g_{1}g)u^{-1}(g)))
\end{array}
}

Note that this identity does not depend on the particular matrix $B^{a}[1]$,
thus, this means that we only need to verify the identity for arbitrary $a,g$ and $g_{1}$:

\ff{
\begin{array}{c}
  \rho(u(g_{1}ga)u^{-1}(ga)))\circ
  (\rho(u(ga)u^{-1}(a))\circ C(g'(a))\circ\rho(u^{-1}(g)))) =\\\\
  =(\rho(u(g_{1}ga)u^{-1}(a))\circ C(g'(a))\circ\rho(u^{-1}(g_{1}g))))
  \circ(\rho(u(g_{1}g)u^{-1}(g))),
\end{array}
}

which is obvious, after the natural simplifications

\ff{
\begin{array}{c}
  \rho(u(g_{1}ga)u^{-1}(a))\circ C(g'(a))\circ\rho(u^{-1}(g)))) =\\\\
  =(\rho(u(g_{1}ga)u^{-1}(a))\circ C(g'(a))\circ\rho(u^{-1}(g))),
\end{array}
}

So, it follows, that for every element $[a]\in G_{0}$ there exist an element
$(A^{a}, a)\in \mathrm{Aut}_G\(G_0\times \(F\otimes V\)\)$.
This means that the homomorphism
\ff{\mathrm{Aut}_G\(G_0\times \(F\otimes V\)\)
\mapr{pr} G_0}
is in fact an epimorphism, and the lemma is proved. $\hfill\blacksquare$

There is an  equivalence between $G$-vector bundles with fiber
$G_0\times \(F\otimes V\) $ over a (compact)
base $X$, where $G$ acts trivially over the base and canonically over the fiber,
and homotopy classes of mappings from $X$ to the space
$B\mathrm{Aut}_G\(G_0\times \(F\otimes V\)\)$.

Lets denote by $\mathrm{Vect}_G(M,\rho)$ the category of $G$-equivariant vector
bundles $\xi=\xi_{0}\otimes V$ with base $M$, where  the action of the group $G$
is quasi-free over the base with finite normal stationary subgroup  $H<G$,
 the group $H$ acts trivially over the bundle $\xi_0$ and $V$ denotes
 the trivial bundle with fiber $V$ and with fiberwise action
of the group $H$ given by an irreducible linear representation $\rho$.
Here we need to require for the representations
\f{\rho_g(h)=\rho(g^{-1}hg)} to be equivalent for every $g\in G$, in the other case,
in view of the commutation rule, this category may be void.

The objects of the category $\mathrm{Vect}_G(M,\rho)$ are ordinary vector
bundles over the space $M$, after tensor product by the fixed bundle $V$.
The action of the group $G$ over this bundles is defined after the tensor
product and the inclusion
$GL(F)\hookrightarrow \mathrm{Aut}_G\(G_0\times \(F\otimes V\)\)$
from lemma 2 ensures that the identities are included.

\begin{corollary}\label{cor1} There is a one-to-one correspondence
\feq{vect}{\mathrm{Vect}_G(M,\rho)\approx [M,BU(F)]}
\end{corollary}$\hfill\blacksquare$

Denote by $\mathrm{Bundle}(X,L)$ the category of principal $L$-bundles over the
base $X$.

\begin{theorem}\label{mono} There is a monomorphism
\feq{categ}{\mathrm{Vect}_G(M,\rho)\longrightarrow
\mathrm{Bundle}(M/G_0,\mathrm{Aut}_G\(G_0\times \(F\otimes V\)\)). }
\end{theorem}

\proof\, By corollary \ref{transitionfunctions}, every element $\xi\in
\mathrm{Vect}_G(M,\rho)$ is defined by transition functions \ff{
\Psi_{\alpha\beta}:\; \(U_\alpha\cap [g_{\alpha\beta}]U_\beta\)\to
\mathrm{Aut}_G\(G_0\times \(F\otimes V\)\)}where, by construction,
when $[g]\neq[g_{\alpha\beta}]$, we have
$U_{\alpha}\cap[g]U_{\beta}=\emptyset$ and if $[g]\neq 1$, then
$U_{\alpha}\cap[g]U_{\alpha}=\emptyset$ and
 $U_{\beta}\cap[g]U_{\beta}=\emptyset$. This means that
 the sets $U_{\alpha}$ and $U_{\beta}$  project homeomorphically to open sets under
 the natural projection $M\to M/G_0$. So, these transition
 functions are well-defined over an atlas of the quotient space $M/G_0$ and
 they form a $G$-bundle with fiber $G_0\times \(F\otimes V\)$
 over  this quotient space.

By the same arguments, it is obvious that every $G$-equivariant map
 \feq{tr1}{h_\alpha:O_\alpha\times \(F\otimes V\)\to O_\alpha\times \(F\otimes V\)}can
be interpreted as a map
\feq{tr2}{h_\alpha:U_\alpha\times\(G_0\times \(F\otimes V\)\)\to U_\alpha\times\(G_0\times \(F\otimes V\)\)}by
means of the homeomorphism \f{O_\alpha\approx U_\alpha\times G_0}, where the set $U_\alpha$
can be thought as an open set of the space $M/G_0$. Equivalently,
\feq{tr3}{h_\alpha:U_\alpha \to \mathrm{Aut}_G\(G_0\times \(F\otimes V\)\)}where
$U_\alpha$ is homeomorphic to an open set of the space $M/G_0$.
Therefore, the map
(\ref{categ}) is well defined.

Conversely, if we start from mappings of the form
(\ref{tr3}) where the sets $U_\alpha$ are open in $M/G_0$, by refining the
atlas, if it is necessary, we can always think that the inverse image of the open sets $U_\alpha$
under the quotient map $M\to M/G_0$ are homeomorphic to the product $U_\alpha\times G_0$
and then obtain mappings of the form (\ref{tr1}). Therefore, the map
(\ref{categ}) is a monomorphism. $\hfill\blacksquare$

Of course, the map (\ref{categ}) its not in general an epimorphism, since, when we define the
category $\mathrm{Vect}_G(M,\rho)$, we are automatically fixing a bundle $M\to M/G_0$, or
equivalently, a homotopy class in $[M/G_0, BG_0]$.

\begin{theorem}\label{clasification} If the space $X$ is compact, then
\feq{union}{\mathrm{Bundle}(X,\mathrm{Aut}_G\(G_0\times \(F\otimes
V\)\))\approx \bigsqcup_{M\in \mathrm{Bundle}(X,G_0)}
\mathrm{Vect}_G(M,\rho).}
\end{theorem}

\proof By theorem \ref{mono}, there is an inclusion \feq{inclusion}{
\bigcup_{M\in \mathrm{Bundle}(X,G_0)}
\mathrm{Vect}_G(M,\rho)\hookrightarrow\mathrm{Bundle}(X,\mathrm{Aut}_G\(G_0\times
\(F\otimes V\)\)).}

Now we will construct an inverse to the map (\ref{inclusion}), so the fact that the last union is disjoint will follow. Let
\ff{ \Psi_{\alpha\beta}:\;
\(U_\alpha\cap U_\beta\)\to \mathrm{Aut}_G\(G_0\times \(F\otimes V\)\)}
be the transition functions of a bundle $\xi\in\mathrm{Bundle}(X,\mathrm{Aut}_G\(G_0\times \(F\otimes V\)\))$.
By lemma \ref{shortexact}, there is a continuous projection
of groups \linebreak $pr:\mathrm{Aut}_G\(G_0\times \(F\otimes V\)\)\to G_0$. So, by composition with
$pr$ we obtain a bundle with the discrete fiber $G_0$, and it is well known that $G_0$ acts
fiberwise and freely over the total space $M$ of this bundle and that $ M/G_0=X$.

 Also, we can assume that we have chosen an atlas such that there is a homeomorphism
 \ff{M\approx \underset{\alpha}{\bigcup}\(U_\alpha\times G_0\)\approx\underset{\alpha}{\bigcup}\(\bigsqcup_{[g]\in G_{0}}[g]U_{\alpha}\)}
 where the intersections are defined by the rule
 \ff{[1]U_\alpha\cap [g_{\alpha\beta}]U_\beta \approx U_\alpha\cap U_\beta}
 where $[g_{\alpha\beta}]=pr\circ \Psi_{\alpha\beta}$.

 On the other hand, we have
 \ff{\xi\approx \underset{\alpha}{\bigcup}\(U_\alpha\times \(G_0\times \(F\otimes V\)\)\)}
 where $U_\alpha\times \(G_0\times \(F\otimes V\)\)$ intersects $U_\beta\times \(G_0\times \(F\otimes V\)\)$
on the points $(x,g,f\otimes v)=(x,\Psi_{\alpha\beta}([g],f\otimes v))=(x,[g_{\alpha\beta}g],A_{\alpha\beta}[g](f\otimes v))$
where
\f{x\in U_\alpha\cap U_\beta}
and, once again,  we are using lemma \ref{shortexact}
for the description of the operators $\Psi_{\alpha\beta}$.

Taking into account the homeomorphism
\ff{U_\alpha\times G_0\approx\bigsqcup_{[g]\in G_{0}}[g]U_{\alpha}}
we can rewrite
\ff{([g]x,f\otimes v)=([gg_{\alpha\beta}]x,A_{\alpha\beta}[g](f\otimes v))}.

 Therefore, the projection
 \ff{\(U_\alpha\times G_0\)\times \(F\otimes V\)\to U_\alpha\times G_0}
 extends to a well-defined and continuous projection
 \ff{\xi\to  M.}It is clear by the preceding formulas,
 that this projection will be $G$-equivariant, if $G$ acts canonically over
  the fibers and in by left translations over its factor $G_0$.
 So, we have $\xi\in \mathrm{Vect}_G(M,\rho)$.

 To end the proof, we make the remark that, by the theory of principal $G_0$-bundles, the construction
 of the space $M$ is up to equivariant homeomorphism. This means that the inverse to
 (\ref{inclusion}) is well defined. $\hfill\blacksquare$

 \begin{corollary}
 If the space  $X$ is compact, then
\feq{union}{[X,B(\mathrm{Aut}_G\(G_0\times \(F\otimes
V\)\))]\approx \bigsqcup_{M\in [X,BG_0]}
[M,BU(F)].}
\end{corollary}
\proof\, This corollary follows from general theory.
It is necessary only to prove,
that $\mathrm{Vect}_G(M,\rho)\not=\emptyset$, for every
$M\in $ $\mathrm{Bundle}(X,BG_0)$, but, for example, we have
\ff{M\times F\otimes V\in \mathrm{Vect}_G(M,\rho).}$\hfill\blacksquare$

\subsection{The case when the subgroup $H<G$ is not normal}

Consider an equivariant vector $G$-bundle $\xi$ over the base $M$
\ff{
\begin{array}{c}
          \xi \\
          \mapd{p} \\
          M. \\
        \end{array}
}Let  $H<G$ be a finite subgroup. Assume
that $M$ is the set of fixed points of the conjugation class of this
subgroup, more accurately
\feq{fixconjugate}{M=\underset{[g]\in G/N(H)}{\bigcup} M^{gHg^{-1}},}
and that there is no more fixed points of the conjugation class of $H$ in
the total space of the bundle $\xi$; here we have denoted by
\f{M^H} the set of fixed points of the action of the subgroup
$H$ over the space $M$, $N(H)$ the normalizer of the group  $H$ in $G$.

In other words,  $H<G$ is the unique, up to conjugation,
 maximal subgroup for the $G$-bundle $\xi$, i.e. $\max\rH_\xi= [H]$.

\begin{lemma}\label{disjointbundles}If the condition (\ref{fixconjugate}) holds,
then the $G$-bundle $\xi$ can be presented as a
disjoint union of pair-wise isomorphic bundles with
quasi-free action over the base. More precisely
\ff{\xi=\underset{[g]\in G/N(H)}{\bigsqcup} \xi_{[g]},}where
\ff{\xi_{[g]}=p^{-1}(M^{gHg^{-1}})}
 is a vector bundle with quasi-free
action of the group $N(gHg^{-1})$ and, for every element $g\in G$
  the mapping
\ff{g:\xi^{[1]}\mapr{}g\xi^{[1]}=\xi^{[g]}}
defines an equivariant isomorphism of this bundles, i.e.
the diagram
 \feq{equivariant}{
\begin{array}{ccc}
  N(H)\times \xi_{[1]} &\mapr{} & \xi_{[1]} \\
  \mapd{s_g\times g} &&\mapd{g}\\
  N(gHg^{-1})\times \xi_{[g]}&\mapr{}& \xi_{[g]}\\
\end{array}}commutes, where
\ff{s_g:N(H)\mapr{}N(gHg^{-1})=gN(H)g^{-1},\quad(g,n)\mapsto gng^{-1}.}
\end{lemma}

\proof\,By the maximality of $H$, the union (\ref{fixconjugate})
is pair-wise disjoint, i.e.
 \ff{M=\underset{[g]\in G/N(H)}{\bigsqcup} M^{gHg^{-1}}}and,
 therefore,
 \ff{\xi=\underset{[g]\in G/N(H)}{\bigsqcup} \xi_{[g]}.}Since
the action of  $G$ is fiberwise, we have $g\cdot \xi_{[1]}=\xi_{[g]}$
 for every $g\in G$. Restricting
 the projection $\xi\mapr{} M$ to the space \f{\xi_{[g]}}, we
 obtain the bundle
\ff{\begin{array}{c}
          \xi_{[g]} \\
          \mapd{p} \\
          M^{gHg^{-1}}. \\
          \end{array}}

 The bundle
 $\xi_{[g]}$ has an action of the normalizer $N(gHg^{-1})$:
 \ff{N(gHg^{-1})\times \xi_{[g]}\longrightarrow \xi_{[g]},} i.e.
 \f{\xi_{[g]}} is a $N(gHg^{-1})$-bundle for every $g\in G$.

 Note that group conjugation
$s_g:N(H)\mapr{} N(gHg^{-1})$ defines an isomorphism
 between these groups that fits into the commutative
 diagram
 \ff{
\begin{array}{ccc}
  N(H)\times \xi_{[1]} &\mapr{} & \xi_{[1]} \\
  \mapd{} &&\mapd{}\\
  N(gHg^{-1})\times \xi_{[g]}&\mapr{}& \xi_{[g]}.\\
\end{array}}i.e. $gng^{-1}\cdot gx=g\cdot nx$. This means that
the bundles \f{\xi_{[1]}} and \f{\xi_{[g]}}
are naturally and equivariantly isomorphic.

Evidently, the mappings on the diagram (\ref{equivariant})
do not depend on the elements $n\in N(H)$,
but they depend on the element $g\in G$.

The action of the group
$N(H)$ over the base $M^H$ reduces to the factor group $N(H)/H$:
\ff{
\begin{array}{ccc}
  N(H)\times  \xi_{[1]}&\mapr{} &\xi_{[1]} \\
  \mapd{} &&\mapd{}\\
  N(H)/H\times M^H&\mapr{}& M^H\\
\end{array}}where, considering the maximality of the group $H$,
the action $N(H)/H\times M^H\mapr{} M^H$ is free and, by hypothesis,
there is no more fixed of the action of the subgroup
$H$ in the total space of the bundle $\xi$, i.e.
$N(H)$  acts quasi-freely over the base and has normal stationary subgroup $H$.$\hfill\blacksquare$

\begin{definition}
If the condition (\ref{fixconjugate}) holds, we will say that the group
$G$ acts quasi-freely over the bundle $\xi$
with (non-normal) stationary subgroup $H$.
\end{definition}

Let $X(\rho)$ be the canonical model for the representation
\f{\rho:H\mapr{} GL(F)} with action of the group $N(H)$. Define
a canonical model  $X(\rho_g)$ for the representation
\ff{\rho_g:gHg^{-1}\mapr{s_{g^{-1}}} H\mapr{\rho} GL(F),}$s_g(n)
=gng^{-1}$. The action of the group $N(gHg^{-1})$ over $X(\rho_g)$ is defined
 using the homomorphism of right $gHg^{-1}$-modules
\ff{u_g:gHg^{-1}\mapr{s_{g^{-1}}}H\mapr{u}N(H)\mapr{s_g} N(gHg^{-1})}by
the formula
\feq{canaction}{
\begin{array}{ll}
\phi([g],g_{1})
=
\bigoplus_k\(\id\otimes\rho_k(u(g_{1}g)u^{-1}(g))\)=
\rho(u(g_{1}g)u^{-1}(g)).
    \end{array}
}

Let
\ff{GX(\rho):=\underset{[g]\in G/N(H)}{\bigsqcup} X(\rho_g)}i.e.
if  $lHl^{-1}=gHg^{-1}$, then the spaces $X(\rho_g)$ and $X(\rho_l)$
coincide.

This notation will be clear after the next lemma.

\begin{lemma}\label{nonnormal}
The group $G$ acts over the space $GX(\rho)$  quasi-freely with
(non-normal) stationary subgroup $H$ and, under this action, the space
$GX(\rho)$ coincides with the orbit of the subspace $X(\rho)$. In particular,
we have the relations
\ff{N(H)\(X(\rho)\)=X(\rho)}and
\ff{(GX(\rho))^{gHg^{-1}}=N(gHg^{-1})/gHg^{-1}.}
\end{lemma}
\proof\,
The action
$G\times GX(\rho)\to GX(\rho)$ is defined
in the following way:
for a fixed $g\in G$ define the mapping
\ff{ g :X(\rho)\mapr{} X(\rho_g) }as
\ff{s_g\times \id :N(H)_0\times F\mapr{}
                       N(gHg^{-1})_0\times F}
(\f{N(H)_0=N(H)/H})
and, if $lHl^{-1}=gHg^{-1}$, then the mapping $l:X(\rho)\mapr{}X(\rho_l)$
is chosen to make the diagramm
\feq{canaction2}{\begin{array}{ccc}
   X(\rho_g)  &\mapr{s_g^{-1}\times \id}& X(\rho)\\
   \|   &                                          & \mapd{l^{-1}g}  \\
  X(\rho_{l}) &\mapr{l^{-1}}                      & X(\rho).
\end{array}}commutative, i.e.
\ff{l=(s_g\times \id)\circ(g^{-1}l)}where
the mapping
\feq{intern}{g^{-1}l:X(\rho)\mapr{}X(\rho)=X(\rho_{g^{-1}l})} is
the canonical left translation by the element $g^{-1}l\in N(H)$.
$\hfill\blacksquare$

\begin{corollary}There is an isomorphism
\feq{isom}{g:\mathrm{Aut}_{N(H)}\(X(\rho)\)\mapr{\approx}
                    \mathrm{Aut}_{N(gHg^{-1})}\(X(\rho_g)\)}that depends
only on the class $[g]\in G/N(H)$.
\end{corollary}

\proof\, We have a diagram  (\ref{equivariant}) for $\xi=GX(\rho)$.
Such a diagram always induces an isomorphism
\ff{\mathrm{Aut}_{N(H)}\(X(\rho)\)\mapr{\approx}
                    \mathrm{Aut}_{N(gHg^{-1})}\(X(\rho_g)\)}by
the rule
\ff{\mathbf{A}\mapsto g\mathbf{A}g^{-1}}and,
if $l\in [g]\in  G/N(H)$ then $l^{-1}g\in N(H)$ commutes with
$\mathbf{A}\in\mathrm{Aut}_{N(H)}\(X(\rho)\)$. Therefore
\ff{g\mathbf{A}g^{-1}=g(g^{-1}l)(l^{-1}g)\mathbf{A}g^{-1}=g(g^{-1}l)\mathbf{A}(l^{-1}g)g^{-1}=l\mathbf{A}l^{-1}.} $\hfill\blacksquare$

\begin{definition}
The space $GX(\rho)$ is called the canonical model for the case
when the subgroup  $H<G$ is not normal.
\end{definition}

\begin{lemma}
\feq{isoaut}{\mathrm{Aut}_{G}\(GX(\rho)\)\approx
          \mathrm{Aut}_{N(H)}\(X(\rho)\)}
\end{lemma}

\proof\,
By definition, an element of the group \f{\mathrm{Aut}_G\(X\)}
is an equivariant mapping $\mathbf{A}^{a}$ such that the pair
$(\mathbf{A}^{a},a)$
defines the commutative diagram
\ff{
\begin{array}{ccc}
   X &\mapr{\mathbf{A}^{a}}& X\\
  \mapd{} &&\mapd{}\\
   G/H&\mapr{a}&G/H, \\
\end{array}
}that
commutes with the canonical action, i.e. the mapping
 \f{a\in \mathrm{Aut}_G(G/H)} satisfies the condition
\ff{a\in \mathrm{Aut}_G(G/H)\approx N(H)/H,\quad a[g]=[ga],\;[g]\in N(H)/H.}
Therefore,
$\mathbf{A}^{a}=(A^{a}[g])_{[g]\in N(H)/H}\in \mathrm{Aut}_{N(H)}(X(\rho))$.

The value of the operators
$(A^{a}[g])_{[g]\in G/H}$
can be calculated in terms of the operator
$A^{a}[1]$ as in lemma 5 (formula \ref{16.d}).
$\hfill\blacksquare$

Denote by
$\widetilde{\mathrm{Vect}}_G(M,\rho)$ the category of vector bundles
with quasi-free action of the group $G$ over the base $M$.

\begin{theorem}\label{redtonormal}
$\widetilde{\mathrm{Vect}}_G(M,\rho)\approx \mathrm{Vect}_{N(H)}(M^H,\rho)$.
\end{theorem}

\proof\,
From lemma \ref{disjointbundles} follows that the bundles
$\xi_{[1]}$ and $\xi_{[g]}$ equivariantly isomorphic and
are given by mappings
\ff{
M^{gHg^{-1}}/N(gHg^{-1})_0\mapr{} B\mathrm{Aut}_{N(gHg^{-1})}\(X(\rho_g)\),
}and
\ff{
M^H/N(H)_0\mapr{} B\mathrm{Aut}_{N(H)}\(X(\rho)\),
}that can be put in the
commutative diagram
\ff{\begin{array}{ccc}
   M^{H}/N(H)_0&\mapr{}& B\mathrm{Aut}_{N(H)}\(X(\rho)\)\\
   \mapd{\bar g}& &\mapd{\bar g}\\
  M^{gHg^{-1}}/N(gHg^{-1})_0&\mapr{}& B\mathrm{Aut}_{N(gHg^{-1})}\(X(\rho_g)\).
\end{array}
}Here, $g:\xi^{H}\mapr{}\xi^{gHg^{-1}}$ is the action over the bundle
 $\xi$. The arrow on the right side is induced by the isomorphism
(\ref{isom}) and does not depend on the element $g\in [g]\in G/N(H)$.
$\hfill\blacksquare$

\section{The computation of the groups $ {}_{n}\Omega_{G}^{v}(\rH_{k},\rH'_{k})$}
Recall that the bordism groups $ {}_{n}\Omega_{G}^{v}(\rH_{k},\rH'_{k})$
are generated by manifolds $M^H$ of fixed points sets of subgroups
$H\in \rH_{k}= \max\{\rH'_{k-1}\}$. Due to the maximality, for subgroups
$H,K\in \rH_{k}= \max\{\rH'_{k-1}\}$, the corresponding manifolds do not
intersect, i.e. $M^K\cap M^H=\emptyset $. So, we can write
$$ {}_{n}\Omega_{G}^{v}(\rH_{k},\rH'_{k})\approx
\bigoplus_{[H]\in \rH_{k}/G} {}_n\Omega_{G}^{v}([H])$$
where $[H]\subset \rH_{k}$ denotes de conjugation class
of the subgroup $H\in \rH_{k}$.

Now, every manifold $M^H$ has free action of the group $N(H)/H$
and is equipped with the structure of its normal bundle with
action of the group $N(H)$.

Lets denote by
\f{{}_{n}\Omega_{N(H)/H}(B\mathrm{Aut}_{N(H)}\(X(\rho)\))},
the group of bordisms of the space \f{B\mathrm{Aut}_{N(H)}\(X(\rho)\)} with
the restriction that the induced map
\ff{\pi_1(pr_*\circ f):\pi_1(X)\longrightarrow \pi_1(BN(H)_0)\approx N(H)_0}
is a monomorphism, where
\f{pr:\mathrm{Aut}_{N(H)}\(X(\rho)\)\mapr{} G_0} is the epimorphism
in the short exact sequence (\ref{suc}).

\begin{theorem} It takes place the isomorphism
\feq{classification}{{}_n\Omega_{G}^{v}([H])\approx
\bigoplus_\rho {}_{n-r}\Omega_{N(H)/H}(B\mathrm{Aut}_{N(H)}\(X(\rho)\))}where
$\rho:H\to U(F),\;r=\dim F$ runs over all the unitary non-equivalent representations
of the group $H$.
\end{theorem}
\proof By theorem \ref{redtonormal}, it is enough to consider a manifold equipped
with a vector bundle $p:\nu_{N(H)}\mapr{} M^H$ with quasi-free action
of the group \f{N(H)} over the base and stationary subgroup $H$.

According to \ref{clasification}, the manifold  $(M^H,\nu_{N(H)})$
is defined by the pair $(M^H,f)$, where $M^H$ is some manifold
with free action of the group \f{N(H)/H} and some
continuous map
\ff{f:M^H/N(H)_0\to B\mathrm{Aut}_{N(H)}\(X(\rho)\),}
which is defined up to homotopy, \f{N(H)_0=N(H)/H},
$\rho:H\to U(F)$ is a linear representation of the group
$H$. In fact, the representation $\rho$ depends on the
connectedness components in the orbit space \f{M^H/N(H)},
but, by simplicity, we will assume that it is fixed.

Lets show that the homomorphism
\ff{\pi_1(pr_*\circ f):\pi_1(X)\mapr{} \pi_1(BN(H)_0)}is
a monomorphism. In order to do so, lets consider the manifold
$M^H$ as the union of its connectedness components:
\ff{M^H=\bigsqcup_sM_s^H.}Then the subgroups
\ff{N(H)_0^s=\{n\in N(H)_0|n(M_s^H)=M_s^H\}}
act freely over the spaces $M_s^H$ and,
therefore
\ff{N(H)_0^s\approx
\pi_1(M_s^H/N(H)_0^s)\longrightarrow \pi_1(BN(H)_0)
\approx N(H)_0} is a monomorphism, but
\f{\pi_1(M_i^H)\approx\pi_1(n(M_n^H))}
for every $n\in N(H)_0$,
\ff{p:p^{-1}(N(H)_0M_i^H)\mapr{} N(H)_0M_i^H}is a vector bundle
with quasi-free action of the group $N(H)$ over the base.

Now, consider a manifold $X$ and a map
\ff{f:X\to B\mathrm{Aut}_{N(H)}\(X(\rho)\)}such that
\ff{\pi_1(pr_*\circ f):\pi_1(X)\mapr{} \pi_1(BN(H)_0)}is a
monomorphism. Then, over the space \f{M=f^*EN(H)_0}, the group $N(H)_0$
acts freely, so this space inherits the smooth structure
from the manifold $X$.

Of course, the corresponding vector bundle with quasi-free action of
the normalizer $N(H)$ is obtained as
$\xi=f^*E\mathrm{Aut}_{N(H)}\(X(\rho)\)$.$\hfill\blacksquare$

\subsection*{Other open problems}
%{\color{red} This part should be shifted to the end of the article}

\begin{itemize}
\item
Compare the index of the signature differential operator on the manifold with proper action with the noncommutative signature of the manifold.

\item Calculate the noncommutative signature on a manifold with proper action
in the terms of fixed points.
\end{itemize}


\begin{thebibliography}{20}
\bibitem{Mishchenko-1998} Luke G., Mishchenko A. S., \textit{Vector Bundles And Their Applications.} Kluwer Academic Publishers Group (Netherlands), 1998, ISBN: 9780792351542
\bibitem{Conner-1964} P. Conner, E. Floyd. \textit{Differentiable periodic maps.} Berlin, Springer-Verlag 1964.
\bibitem{Atiyah-1967} Atiyah M.F., \textit{K-theory.} Benjamin, New York, (1967).
\bibitem{Atiyah-Segal-1965} Atiyah, M.F., Segal, G.B., \textit{Equivariant K-theory Lecture Notes} Oxford, (1965).

\bibitem{BaCoHi-94e}
P.Baum, A.Connes, and N.Higson.
\newblock Classifying space for proper actions and $K$-theory of group
  $C^*$-algebras.
\newblock {\em Contemp. Math.}, 167:241--291, 1994.

\bibitem{Levine-2008} Levine M.,  Serpe C.,\textit{On a spectral sequence for equivariant K-theory}
    K-Theory (2008) 38, No.2, pp. 177–222, arXiv:math/0511394v3 [math.KT] 19 Nov 2005.

\bibitem{Paw-02e}
K.~Pawalowski. \textit{Manifolds as the fixed point sets of smooth compact lie group
  actions.} Current Trends in Transformation Groups, K-Monographs in
  Mathematics 7, pages 79--104. Kluwer Academic Publishers.

\bibitem{Ill-00e}
S.~Illman.
\newblock Existence and uniqueness of equivariant triangulations of smooth
  proper gmanifolds with some applications to equivariant whitehead torsion.
\newblock {\em J. Reine Angew. Math.}, 524:129--183, 2000.



\bibitem{Kor-05en}
T.~Korppi.
\newblock Equivariant triangulations of differentiable and real-analytic
  manifolds with a properly discontinuous action.
\newblock In {\em Annales Academi acientiarum fennic matematica
  dissertationes,}, number 141. Suomalainen Tiedeakatemia, XVC05/4352.
  b20453085., Helsinki, 2005.



\bibitem{MishQuit}Mishchenko, A.S. and Morales Melendez, Q.
\textit{Description of the vector $G$-bundles over $G$-spaces
with quasi-free proper action of discrete group $G$.}
arXiv:0901.3308v1  [math.AT]  21 Jan 2009

\bibitem{Quit}Morales Melendez, Q.
\textit{Description of G-bundles over G-spaces with
quasi-free proper action of discrete group II.}
arXiv:0912.5047v1  [math.KT]  27 Dec 2009

\bibitem{MishQuit-1}
Mishchenko, A.S. and Morales Melendez, Q.
\textit{Description of the  $G$-bundles over $G$-spaces with
quasi-free proper action of discrete group $G$.} (in Russian)
Deposited at VINITI, no. 716-V2009, 24.11.2009, 15 pp.

\bibitem{[2]}{Morales Melendez, Q. \textit{Description of the
$G$-bundles over $G$-spaces with
quasi-free proper action of discrete group $G$.
II.}  (in Russian)
Deposited at VINITI, no. 717-V2009, 24.11.2009, 11 pp.}

\bibitem{[3]}{Morales Melendez, Q. \textit{Bordisms of Manifolds
with Proper Action of a Discrete Group.}  (in Russian)
Deposited at VINITI, no. 718-V2009, 24.11.2009, 13 pp.}



\bibitem{[1]} {M.K.Morales \textit{Bordisms of Manifolds with Proper
Action of a Discrete Group} (in Russian)
Vestnik Moskovskogo Universiteta. Seria 1. Matematika. Mekhanika. 2010. no. 2. p. 92--94.
Englisj translation:
Moscow University Mathematics Bulletin,
Volume 65, Number 2, 92--94, (2010), DOI: 10.3103/S0027132210020105 }

\bibitem{[4]}
Mishchenko, A.S,, Morales Melendez, Q., \textit{
Description of the
$G$-bundles over $G$-spaces with
quasi-free proper action of discrete group $G$.}(in Russian)
Proceedings of the seminar by P.K.Rashevski on vector and 
tensor analysis. With applications to geometry, mechanics and physics. Vyp.27,
p. 150--172 (2011).






\end{thebibliography}
\end{document}